\documentclass[11pt]{article}
\usepackage{amssymb,amsfonts,amsmath,latexsym,epsf,tikz,url}
\usepackage{graphicx}
\usepackage[font={small}]{caption}
\usepackage{cite}

\newtheorem{theorem}{Theorem}[section]

\newtheorem{corollary}[theorem]{Corollary}
\newtheorem{lemma}[theorem]{Lemma}

\newcommand{\proof}{\noindent{\bf Proof.\ }}
\newcommand{\qed}{\hfill $\square$\medskip}

\begin{document}

\title{\textbf{Diminished Sombor index and its relationship with topological indices}}
\author{Fateme Movahedi \footnote{Corresponding author \, E-mail: f.movahedi@gu.ac.ir}}

\maketitle

\begin{center}
Department of Mathematics, Faculty of Sciences, Golestan University, Gorgan, Iran.
\end{center}
\maketitle

\begin{abstract}
In this paper, we investigate the Diminished Sombor index (DSO), a recently introduced degree-based topological index for a simple graph $G$, defined as
\[
DSO(G) = \sum_{uv \in E} \frac{\sqrt{d_u^2+d_v^2}}{d_u+d_v},
\]
where $d_u$ denotes the degree of a vertex $u \in V$. We establish several sharp bounds for this index in terms of classical topological indices such as the Zagreb, Albertson, Harmonic, Randi\'c, and geometric-arithmetic indices. The relationships and inequalities between DSO and these indices are analyzed thoroughly, with characterizations of extremal graphs achieving equality conditions.
\end{abstract}

\medskip
\noindent{\bf AMS Subj.\ Class.:} 05C09, 05C92, 05C90.

\section{Introduction}
Let $G$ be a simple graph with vertex set $V$ and edge set $E$, where $|V|$ and $|E|$ represent the number of vertices (order) and edges (size) of the graph, respectively. The degree of a vertex $u$, denoted by $d_u$, refers to the number of vertices adjacent to $u$. The maximum and minimum degrees of the graph are denoted by $\Delta$ and $\delta$, respectively. An edge connecting two adjacent vertices $u$ and $v$ is denoted as $uv \in E$. The complement of $G$, denoted by $\bar{G}$, is a graph with the same vertex set $V$, in which two vertices are adjacent if and only if they are not adjacent in $G$. Any additional graph-theoretic terminology used in this paper but not defined here, can be found in \cite{1}.

In mathematical chemistry, topological indices, numerical descriptors derived from molecular graphs, have emerged as indispensable tools. These indices quantitatively capture structural features of molecules and play a key role in uncovering correlations between molecular structure and various physicochemical properties. Their predictive power, which eliminates the need for experimental synthesis, makes them highly effective for virtual screening in the search for potential drug candidates or materials with desired characteristics \cite{top}. For a detailed exploration of degree-based topological indices, readers are referred to the comprehensive survey in \cite{2}.\\

The first Zagreb index ($M_1(G)$) of a graph $G$ formulated by Gutman and Trinajsti\'c \cite{3} as follows 
$$M_1(G)=\sum_{uv \in E} \left(d(u)+d(v)\right)=\sum_{u \in V}d(u)^2.$$ 
In 1975, Milan Randi\'c \cite{Rand} introduced the Randi\'c index which is defined as follows
$$R(G)=\sum_{uv \in E}\frac{1}{\sqrt{d_u d_v}}.$$
In \cite{Alb}, the Albertson index which is sometimes referred to as the third Zagreb index \cite{M3} is defined as 
$$Alb(G)=\sum_{uv \in E}|d_u-d_v|.$$ 
The geometric-arithmetic index (GA index), introduced by Vuki\^cevi\^c and Furtula in \cite{GA}, is defined as
$$GA(G)=\sum_{uv \in E}\frac{2\sqrt{d_u d_v}}{d_u+d_v}.$$ 
The Harmonic index is defined in \cite{Har} as 
$$H(G)=\sum_{uv \in E}\frac{2}{d_u+d_v},$$
and in \cite{ISI} the inverse sum indeg index is defined as 
$$ISI(G)=\sum_{uv \in E}\frac{d_u d_v}{d_u+d_v}.$$ 
Kulli in \cite{GAF} introduced the geometric-arithmetic F-index of a graph $G$ as 
$$GAF(G)=\sum_{uv \in E}\frac{2d_ud_v}{d_u^2+d_v^2}.$$
The sum-connectivity index ($\chi(G)$) was established in \cite{SCI} and can be expressed as
$$\chi(G)=\sum_{uv\in E}\frac{1}{\sqrt{d_u+d_v}}.$$
A generalization of the $\chi(G)$ index, known as the general sum-connectivity index, was introduced in \cite{GSCI}. It is formally defined as
$$\chi_{\alpha}(G)=\sum_{uv \in E}\left(d_u+d_v\right)^{\alpha},$$
where $\alpha \in \mathbb{R}$ is an arbitrary number. \\
The forgotten index ($F(G)$) \cite{Forg} and the multiplicative forgotten index ($\Pi_F(G)$) \cite{Forgp} are defined as 
$$F(G)=\sum_{uv \in E}(d_u^2+d_v^2), ~~~~~~~~~~~~~~~\Pi_F(G)=\prod_{uv \in E}(d_u^2+d_v^2).$$
The sum-connectivity F-index, proposed by Kulli in \cite{SCF}, is defined as
$$SF(G)=\sum_{uv \in E}\frac{1}{\sqrt{d_u^2+d_v^2}}.$$
In 2021, Kulli introduced the first Banhatti–Sombor index for a connected graph $G$ in \cite{BSI}, defined as:
$$BSO(G)=\sum_{uv \in E}\sqrt{\frac{1}{d_u^2}+\frac{1}{d_v^2}}.$$
The Adriatic indices were introduced and studied over a series of publications several years ago \cite{SDD}. Among them, only a limited number have demonstrated potential utility in modeling physico-chemical properties of molecular structures. One such index is the so-called symmetric division degree index ($SDD(G)$), which is defined as
$$SDD(G)=\sum_{uv \in E}\frac{d_u^2+d_v^2}{d_u d_v}.$$
The Sombor index, a prominent topological index in graph theory, is proposed in \cite{11} and is defined as $SO(G)=\sum_{uv \in E}\sqrt{d_u^2+d_v^2}$, has garnered significant attention for its applications in quantitative structure-property/activity relationships (QSPR/QSAR) studies \cite{11, 12, 13, 14, 15, 16,166,1666}. Recently, a variant of the Sombor index was introduced in \cite{X} and is called the Diminished Sombor index in \cite{FaGut}. This index is defined as follows
$$DSO(G)=\sum_{uv \in E}\frac{\sqrt{d_u^2+d_v^2}}{d_u+d_v}.$$
Movahedi et al. \cite{FaGut} derived bounds for the Diminished Sombor (DSO) index, identified extremal graphs, and established Nordhaus-Gaddum-type inequalities. They also presented numerical analyses regarding the structure-dependence of the DSO index and its potential applications in chemistry. In \cite{Fateme3}, the tricyclic graph of a specified order that attains the maximum DSO is identified, and its distinctive structural characteristics are examined.

In this paper, we aim to advance the study of the DSO index by deriving new bounds in terms of fundamental graph parameters and exploring its connections with several classical topological indices. Our objective is to offer a deeper understanding of the Diminished Sombor index and to encourage its further investigation and application in the field of chemical graph theory.

\section{Preliminaries} 
In this section, we present several known inequalities and results that will be used in the proofs of our main theorems.

\begin{lemma}{\rm \cite{Mit}}\label{lemma5} 
Let $a=(a_i)_{i=1}^n$ and $b=(b_i)_{i=1}^n$  be the real numbers sequences where $a_i\geq 0$ and $b_i>0$. Then for any $t\leq 0$ or $t\geq 1$, we have
$$\Big(\sum_{i=1}^na_i b_i\Big)^t\leq \Big(\sum_{i=1}^n a_i b_i^t\Big)\Big(\sum_{i=1}^n a_i\Big)^{t-1}.$$
Equality holds if and only if $t=0$ or $t=1$ or $b_1=\cdots=b_n$ or $a_1=\cdots=a_r$ and $b_{r+1}=\cdots=b_n$ for $1\leq r \leq n-1$.
\end{lemma}

\begin{lemma}{\rm \cite{Radon}}\label{lemma6}
Let $x=(x_i)_{i=1}^{n}$, and $a=(a_i)_{i=1}^n$ be two positive real number sequences. For any $r > 0$,
$$\frac{\left(\sum_{i=1}^n x_i\right)^{r+1}}{\left(\sum_{i=1}^n a_i\right)^r}\leq \sum_{i=1}^{n}\frac{x_i^{r+1}}{a_i^r}.$$
Equality holds if and only if either $r=0$ or $\frac{x_1}{a_1}=\cdots=\frac{x_n}{a_n}$.
\end{lemma}

\begin{lemma}{\rm \cite{24}}\label{lemma4} 
Let $G$ be a connected graph with $n$ vertices and $m$ edges. Then
$$M_1(G)\leq m(m+1),$$
with equality for $n>3$ if and only if $G\simeq K_3$ or $G\simeq S_n$. 
\end{lemma}

\begin{lemma}{\rm \cite{Mar}}\label{lemma1} 
If $a_i, b_i \geq 0$ and $Ab_i \leq a_i \leq Bb_i$ for $1\leq i \leq n$, then
$$\Big(\sum_{i=1}^n a_i^2 \Big)\Big(\sum_{i=1}^n b_i^2 \Big)\leq \frac{(A+B)^2}{4AB}\Big(\sum_{i=1}^n a_ib_i \Big)^2.$$
The equality holds if and only if $A=B$ and $a_i=Ab_i$.
\end{lemma}

\begin{lemma}{\rm \cite{Ha}}\label{lemma3} 
Let $G$ be a simple graph with $m$ edges. Then
$$H(G)\geq \frac{2m^2}{M_1(G)}.$$
\end{lemma}

\begin{lemma}{\rm \cite{SS}}\label{lemma7} 
Let $\mathbf{a} = (a_i)_{i=1}^n$ and $\mathbf{b} = (b_i)_{i=1}^n$ be two nonnegative sequences arranged in decreasing order, such that $a_1 \neq 0$ and $b_1 \neq 0$, and let $\mathbf{w} = (w_i)_{i=1}^n$ be a nonnegative sequence. Then the following inequality holds
\[
\left( \sum_{i=1}^n w_i a_i^2 \right) \left( \sum_{i=1}^n w_i b_i^2 \right)
\leq 
\max \left( 
b_1 \sum_{i=1}^n w_i a_i,\,
a_1 \sum_{i=1}^n w_i b_i
\right)
\sum_{i=1}^n w_i a_i b_i.
\]
\end{lemma}

\begin{lemma}{\rm \cite{SSS}}\label{lemma8} 
Let $a=(a_i)_{i=1}^n$, be positive real number sequence. Then
$$\left( \sum_{i=1}^n \sqrt{a_i} \right)^2\geq \sum_{i=1}^n a_i +n(n-1)\left( \prod_{i=1}^n a_i \right)^{\frac{1}{n}}.$$
\end{lemma}

\section{Main Results}
In this section, we obtain several bounds for the Diminished Sombor index in terms of some graph parameters.

\begin{theorem}\label{theorem2}
For a simple graph $G$ of order $n\geq 2$ with the maximum degree $\Delta$, 
$$\frac{\sqrt{2}}{4\Delta} Alb(G)+\frac{1}{2}GA(G)\leq DSO(G) \leq \frac{1}{2\delta}Alb(G)+\frac{\sqrt{2}}{2}GA(G).$$
Equality on the left-hand side holds if and only if $G\simeq \overline{K_n}$ and in the right-hand side holds if and only if $G$ consists of components,
each of which is a regular graph (not necessarily of equal degree).
\end{theorem}

\proof
Since for any two real numbers $a,b\geq 0$, we have
\begin{equation}\label{1}
\frac{\sqrt{2}}{2}\left(\sqrt{a}+\sqrt{b}\right)\leq \sqrt{a+b},
\end{equation}
and the equality holds if and only if $a=b$.

Using the inequality (\ref{1}) and by setting $a=\left(d_u-d_v\right)^2$ and $b=2d_u d_v$, we get
\begin{align*}
\frac{\sqrt{2}}{2}\Big(\sqrt{\left(d_u-d_v\right)^2}+\sqrt{2d_ud_v}\Big)&\leq \sqrt{d^2_u+d^2_v}\\
\Longrightarrow \frac{\sqrt{2}}{2}\Big(\frac{|d_u-d_v|}{d_u+d_v}+\sqrt{2}\frac{\sqrt{d_ud_v}}{d_u+d_v}\Big)&\leq\frac{\sqrt{d^2_u+d^2_v}}{d_u+d_v}
\end{align*}
Since $2\delta \leq d_u+d_v\leq 2\Delta$ for any $u, v \in V$, we obtain
\begin{align*}
DSO(G)&=\sum_{uv \in E}\frac{\sqrt{d_u^2+d_v^2}}{d_u+d_v}\\
&\geq \frac{\sqrt{2}}{2}\sum_{uv\in E}\frac{|d_u-d_v|}{d_u+d_v}+\frac{1}{2}\sum_{uv\in E}\frac{2\sqrt{d_ud_v}}{d_u+d_v}\\
&\geq \frac{1}{2\Delta} \Big(\frac{\sqrt{2}}{2}\sum_{uv\in E}|d_u-d_v|\Big)+\frac{1}{2}\sum_{uv\in E}\frac{2\sqrt{d_ud_v}}{d_u+d_v}\\
&=\frac{\sqrt{2}}{4\Delta} Alb(G)+\frac{1}{2}GA(G).
\end{align*}
Equality in the above inequalities holds if and only if $G\simeq \bar{K_n}$. \\
For proving the upper bound, we use this fact that $\sqrt{a+b}\leq \sqrt{a}+\sqrt{b}$ and the equality holds if and only if $a=0$ or $b=0$. Therefore, by substituting $a=(d_u-d_v)^2$ and $b=2d_ud_v$, we obtain
$$\sqrt{d^2_u+d^2_v}\leq \sqrt{(d_u-d_v)^2}+\sqrt{2d_ud_v}=|d_u-d_v|+\sqrt{2}\sqrt{d_u d_v}.$$
Since $2\delta\leq d_u+d_v\leq 2\Delta$, we have
\begin{align*}
DSO(G)&=\sum_{uv \in E}\frac{\sqrt{d_u^2+d_v^2}}{d_u+d_v}\\
&\leq \sum_{uv \in E}\frac{|d_u-d_v|}{d_u+d_v}+\sqrt{2}\sum_{uv \in E}\frac{\sqrt{d_u d_v}}{d_u+d_v}\\
&\leq \frac{1}{2\delta}\sum_{uv \in E}|d_u-d_v|+\frac{\sqrt{2}}{2}\sum_{uv \in E}\frac{2\sqrt{d_u d_v}}{d_u +d_v}\\
&=\frac{1}{2\delta}Alb(G)+\frac{\sqrt{2}}{2}GA(G).
\end{align*}
Equality holds if and only if $d_u=d_v$ for any $uv \in E$.
\qed

\begin{theorem}\label{theorem3}
Let $G$ be a simple graph of size $m$. Then
$$DSO(G) \leq \frac{\sqrt{2}}{2}(Alb(G)+m).$$
Equality holds  if and only if $G$ consists of components, each of which is a regular graph (not necessarily of equal degree).
\end{theorem}

\proof
By substituting $a=\frac{1}{2}(d_u-d_v)^2$ and $b=\frac{1}{2}(d_u+d_v)^2$ in the inequality $\sqrt{a+b}\leq \sqrt{a}+\sqrt{b}$, we have
\begin{align*}
\sqrt{\frac{1}{2}(d_u-d_v)^2+\frac{1}{2}(d_u+d_v)^2}&\leq \frac{\sqrt{2}}{2}\Big(\sqrt{(d_u-d_v)^2}+\sqrt{(d_u+d_v)^2}\Big)\\
\Longrightarrow ~~~~~~\quad \quad \quad \frac{\sqrt{d_u^2+d_v^2}}{d_u+d_v}&\leq \frac{\sqrt{2}}{2}\Big(\frac{|d_u-d_v|}{d_u+d_v}+\frac{d_u+d_v}{d_u+d_v} \Big).
\end{align*}
Therefore, we get
\begin{align*}
DSO(G)&=\sum_{uv \in E}\frac{\sqrt{d_u^2+d_v^2}}{d_u+d_v}\leq \frac{\sqrt{2}}{2}\sum_{uv \in E}\Big(\frac{|d_u-d_v|}{d_u+d_v}+\frac{d_u+d_v}{d_u+d_v} \Big)\\
&\leq \frac{\sqrt{2}}{2}\sum_{uv \in E} |d_u-d_v|+\frac{\sqrt{2}}{2}m=\frac{\sqrt{2}}{2}\Big(Alb(G)+m\Big).
\end{align*}
Equality holds if and only if $d_u=d_v$ for any $uv \in E$.
\qed

\begin{corollary}\label{cor4}
Let $G$ be a simple connected graph of size $m$ with the maximum degree $\Delta$. Then 
$$DSO(G)\leq \frac{\sqrt{2}}{2}m(\Delta-\delta+1),$$
with equality holds if and only if $G$ is regular graph.
\end{corollary}
\proof
Since $Alb(G) \leq m(\Delta-\delta)$, using Theorem \ref{theorem3} we have
\begin{align*}
DSO(G) &\leq \frac{\sqrt{2}}{2}(Alb(G)+m)\\
&\leq \frac{\sqrt{2}}{2}(m(\Delta-\delta)+m)\\
&=\frac{\sqrt{2}}{2}m(\Delta-\delta+1).
\end{align*}
The equality holds if and only if $G$ is regular graph.
\qed

\begin{theorem}\label{theorem5}
For a simple graph $G$, 
$$DSO(G)\leq \frac{\sqrt{2}}{2}\left(\sqrt{H(G)\big(M_1(G)-H(G)\big)}\right).$$
Equality holds if and only if sum of the squares of the degrees of adjacent vertices is a contrast.
\end{theorem}
\proof
By substituting $a_i=\frac{1}{d_u+d_v}$ and $b_i=\sqrt{d^2_u+d^2_v}$ into Lemma \ref{lemma5}, we have
\begin{align*}
\big(DSO(G)\big)^2&=\Big(\sum_{uv \in E}\frac{\sqrt{d_u^2+d_v^2}}{d_u+d_v}\Big)^2\\
&\leq \Big(\sum_{uv \in E}\frac{d^2_u+d^2_v}{d_u+d_v}\Big)\Big(\sum_{uv \in E}\frac{1}{d_u+d_v} \Big)\\
&\leq \Big(\sum_{uv \in E}\frac{\left(d_u+d_v\right)^2}{d_u+d_v}-\sum_{uv \in E}\frac{2}{d_u+d_v} \Big)\Big(\frac{1}{2}\sum_{uv \in E}\frac{2}{d_u+d_v}\Big)\\
&=\frac{1}{2}\left(M_1(G)-H(G)\right)H(G).
\end{align*}
Consequently, we have
$$DSO(G)\leq \frac{\sqrt{2}}{2}\left(\sqrt{H(G)\big(M_1(G)-H(G)\big)}\right).$$
Equality holds if and only if $\sqrt{d^2_u+d^2_v}$ is a constant for any $uv \in E$. This implies that the equality holds if and only if $d^2_u+d^2_v$ is a constant for any two adjacent vertices. 
\qed


\begin{theorem}\label{theorem7}
Let $G$ be a simple graph of size $m$. Then
$$DSO(G)\leq \frac{m}{m+1}\sqrt{(m+1)^2-2}.$$
\end{theorem}
\proof
We consider $f(x)=\sqrt{x(\alpha-x)}$ and thus, $f'(x)=\frac{\alpha-2x}{2\sqrt{x\alpha-x^2}}>0$ for $\frac{\alpha}{2}\leq x \leq \alpha$ and $f'(x)<0$ for $x\leq \frac{\alpha}{2}$. Since $H(G)\leq \frac{M_1(G)}{2}$ and using Lemma \ref{lemma3} we have $H(G)\geq \frac{2m^2}{M_1(G)}$, thus $f(H(G))\leq f\left(\frac{2m^2}{M_1(G)}\right)$. Therefore using Theorem \ref{theorem5}, we have
\begin{align*}
DSO(G)&\leq \frac{\sqrt{2}}{2}\left(\sqrt{H(G)\big(M_1(G)-H(G)\big)}\right)\\
&\leq \frac{\sqrt{2}}{2}\left(\sqrt{\frac{2m^2}{M_1(G)}\Big(M_1(G)-\frac{2m^2}{M_1(G)}\Big)}\right)\\
&=m\sqrt{\frac{M_1(G)}{M_1(G)}-\frac{2m^2}{M_1(G)^2}}\\
&=m\sqrt{1-\frac{2m^2}{M_1(G)^2}}.
\end{align*}
On the other hand, we suppose that $g(x)=\sqrt{1-\frac{2m^2}{x^2}}$ and $g'(x)=\frac{m^2x}{\sqrt{1-\frac{2m^2}{M_1(G)^2}}}>0$. Therefore, using Lemma \ref{lemma4}, we get
$$DSO(G)\leq m\sqrt{1-\frac{2m^2}{M_1(G)^2}}=\frac{m}{m+1}\sqrt{(m+1)^2-2}.$$
\qed

\begin{theorem}\label{theorem11}
Let $G$ be a graph of size $m$ with the maximum degree $\Delta$ and the minimum degree $\delta$. Then
$$DSO(G)\geq \sqrt{\alpha m\left(m-GAF(G)\right)},$$
where $\alpha=\frac{2\sqrt{2}(\Delta+\delta)(\sqrt{\Delta^2+\delta^2})}{\left(\sqrt{2(\Delta^2+\delta^2)}+(\Delta+\delta)\right)^2}$. Equality holds if and only if $G$ consists of components, each of which is a regular graph (not necessarily of equal degree).
\end{theorem}
\proof
In the first, we prove the following inequality for any $u, v \in V$
\begin{equation}\label{5}
\frac{1}{\sqrt{2}}\leq \frac{\sqrt{d_u^2+d_v^2}}{d_u+d_v}\leq \frac{\sqrt{\Delta^2+\delta^2}}{\Delta+\delta}.
\end{equation}
Let $f(x)=\frac{\sqrt{1+x^2}}{1+x}$. Since $f'(x)=\frac{(x-1)}{(1+x)^2\sqrt{1+x^2}}>0$ for any $x\geq 1$, thus $f(x)$ is an increasing function on $x\geq 1$. We know that $0<\delta\leq d_u\leq \Delta$ for any $u \in V$, and consequently, $\frac{\delta}{\Delta}\leq \frac{d_u}{d_v}\leq \frac{\Delta}{\delta}$. Since $f(x)$ is the increasing function, we have $f\left(\frac{d_u}{d_v}\right)\leq f\left(\frac{\Delta}{\delta}\right)$. Therefore, 
$$\frac{\sqrt{d_u^2+d_v^2}}{d_u+d_v}=\frac{\sqrt{1+\frac{d^2_u}{d^2_v}}}{1+\frac{d_u}{d_v}}\leq \frac{\sqrt{1+\frac{\Delta^2}{\delta^2}}}{1+\frac{\Delta}{\delta}}=\frac{\sqrt{\Delta^2+\delta^2}}{\Delta+\delta}.$$
On the other hand, since $\left(d_u-d_v\right)^2\geq 0$, we have $\frac{1}{\sqrt{2}}\leq \frac{\sqrt{d_u^2+d_v^2}}{d_u+d_v}$. \\
By substituting $a_i=\frac{\sqrt{d_u^2+d_v^2}}{d_u+d_v}$, $b_i=1$, $A=\frac{1}{\sqrt{2}}$, and $B=\frac{\sqrt{\Delta^2+\delta^2}}{\Delta+\delta}$ into Lemma \ref{lemma1}, we get
 \begin{align*}
\Big(\sum_{uv \in E} \frac{d_u^2+d_v^2}{(d_u+d_v)^2} \Big)\Big(\sum_{uv \in E} 1 \Big)&\leq \frac{(\frac{1}{\sqrt{2}}+\frac{\sqrt{\Delta^2+\delta^2}}{\Delta+\delta})^2}{\frac{4}{\sqrt{2}}\times\frac{\sqrt{\Delta^2+\delta^2}}{\Delta+\delta}}\Big(\sum_{uv \in E} \frac{\sqrt{d_u^2+d_v^2}}{d_u+d_v} \Big)^2\\
\Longrightarrow ~~~~~m\Big(\sum_{uv\in E}1-\sum_{uv\in E}\frac{2d_ud_v}{(d_u+d_v)^2} \Big)&\leq \frac{\big(\sqrt{2(\Delta^2+\delta^2)}+\Delta+\delta\big)^2}{4\sqrt{2}(\Delta+\delta)(\sqrt{\Delta^2+\delta^2})}\big(DSO(G)\big)^2.
 \end{align*}
 Since $GAF(G)=\sum_{uv \in E}\frac{2d_ud_v}{d^2_u+d^2_v}\geq \sum_{uv \in E}\frac{2d_ud_v}{(d_u+d_v)^2}$, we have
 $$m\left(m-GAF(G) \right)\frac{4\sqrt{2}(\Delta+\delta)(\sqrt{\Delta^2+\delta^2})}{\big(\sqrt{2(\Delta^2+\delta^2)}+\Delta+\delta\big)^2}\leq \left(DSO(G)\right)^2.$$
 With considering $\alpha=\frac{4\sqrt{2}(\Delta+\delta)(\sqrt{\Delta^2+\delta^2})}{\left(\sqrt{2(\Delta^2+\delta^2)}+(\Delta+\delta)\right)^2}$,
$$DSO(G)\geq \sqrt{\alpha m\left(m-GAF(G)\right)}.$$
Equality holds if and only if $\frac{1}{\sqrt{2}}= \frac{\sqrt{\Delta^2+\delta^2}}{\Delta+\delta}=\frac{\sqrt{d_u^2+d_v^2}}{d_u+d_v}$ for any $uv \in E$, that is, for any $uv \in E$, $\delta=d_u=d_v=\Delta$. 
\qed

\begin{theorem}\label{theorem12}
Let $G$ be a simple connected graph with the maximum degree $\Delta$ and the minimum degree $\delta$. Then
$$DSO(G)\geq \frac{M_1(G)-2ISI(G)+\delta\Delta H(G)}{\sqrt{2}(\Delta+\delta)}.$$
Equality holds if and only if $G$ is regular graph. 
\end{theorem}
\proof
For any $u, v \in V$, we have $2\delta\leq \sqrt{d_u^2+d_v^2}\leq 2\Delta$. Hence 
$$\left(\sqrt{d_u^2+d_v^2}-2\delta \right)\left(2\Delta-\sqrt{d_u^2+d_v^2} \right)\geq 0.$$ 
Therefore,
$$\sqrt{2}(\Delta+\delta)\sqrt{d_u^2+d_v^2}\geq d^2_u+d^2_v+2\delta\Delta.$$
Using the definition of the DSO index, we get
\begin{align*}
DSO(G)&=\sum_{uv \in E}\frac{\sqrt{d_u^2+d_v^2}}{d_u+d_v}\\
&\geq \frac{1}{\sqrt{2}(\Delta+\delta)}\Big(\sum_{uv \in E}\frac{d^2_u+d^2_v+2\delta\Delta}{d_u+d_v} \Big)\\
&=\frac{1}{\sqrt{2}(\Delta+\delta)}\Big(\sum_{uv \in E}\frac{(d_u+d_v)^2-2d_ud_v}{d_u+d_v}+\delta\Delta\sum_{uv \in E}\frac{2}{d_u+d_v} \Big)\\
&=\frac{1}{\sqrt{2}(\Delta+\delta)}\Big(\sum_{uv \in E}(d_u+d_v)-2\sum_{uv \in E}\frac{d_ud_v}{d_u+d_v}+\delta\Delta\sum_{uv \in E}\frac{2}{d_u+d_v} \Big)\\
&=\frac{1}{\sqrt{2}(\Delta+\delta)}\Big(M_1(G)-2ISI(G)+\delta\Delta H(G)\Big).
\end{align*}
Equality holds if and only if $d_u=d_v=\delta=\Delta$, that is, $G$ is regular graph. 
\qed

\begin{theorem}\label{theorem13}
For a simple $G$, 
$$\frac{(\chi(G))^2}{SF(G)}\leq DSO(G) \leq \sqrt{F(G)\chi_{-2}(G)}$$
Equality holds if and only if $G$ consists of components, each of which is a regular graph (not necessarily of equal degree).
\end{theorem}
\proof
Using the Cauchy-Schwarz inequality, we get
\begin{align*}
\left(DSO(G)\right)^2&=\Big(\sum_{uv\in E}\frac{\sqrt{d_u^2+d_v^2}}{d_u+d_v}\Big)^2\\
&\leq \Big(\sum_{uv \in E}\left(\sqrt{d_u^2+d_v^2}\right)^2\Big)\Big(\sum_{uv \in E}\frac{1}{\left(d_u+d_v\right)^2}\Big)\\
&=\Big(\sum_{uv \in E}\left(d_u^2+d_v^2\right)\Big)\Big(\sum_{uv \in E}\left(d_u+d_v\right)^{-2}\Big)\\
&=F(G)\chi_{-2}(G).
\end{align*}
Hence 
$$DSO(G) \leq \sqrt{F(G)\chi_{-2}(G)}.$$
For the lower bound, we use Lemma \ref{lemma6} in which $x_i=\frac{1}{\sqrt{d_u+d_v}}$, $a_i=\frac{1}{\sqrt{d_u^2+d_v^2}}$, and $r=1$. Therefore, we get
\begin{align*}
DSO(G)&=\sum_{uv\in E}\frac{\sqrt{d_u^2+d_v^2}}{d_u+d_v}=\sum_{uv\in E}\frac{\frac{1}{\left(\sqrt{d_u+d_v}\right)^2}}{\frac{1}{\sqrt{d_u^2+d_v^2}}}\\
&\geq \frac{\Big(\sum_{uv \in E}\frac{1}{\sqrt{d_u+d_v}}\Big)^2}{\sum_{uv \in E}\frac{1}{\sqrt{d_u^2+d_v^2}}}=\frac{\left(\chi(G)\right)^2}{SF(G)}.
\end{align*}
Equality holds if and only if $d_u=d_v$ for any edge $uv \in E$. 
\qed
\begin{theorem}\label{theorem14}
Let $G$ be a connected graph. Then
$$DSO(G)\geq \frac{\left(R(G)\right)^2}{BSO(G)}.$$
\end{theorem}
\proof
By substituting $t=2$, $a_i=\frac{\sqrt{d^2_u+d^2_v}}{d_u+d_v}$ and $b_i=\sqrt{\frac{d_u+d_v}{d_ud_v}}$ into Lemma \ref{lemma5}, we have
\begin{align*}
BSO(G)&=\sum_{uv \in E}\frac{\sqrt{d^2_u+d^2_v}}{d_ud_v}\\
&=\sum_{uv \in E}\left(\frac{\sqrt{d^2_u+d^2_v}}{d_u+d_v}\right)\left(\sqrt{\frac{d_u+d_v}{d_ud_v}}\right)^2\\
&\geq \frac{\left(\sum_{uv \in E}\left(\frac{\sqrt{d^2_u+d^2_v}}{d_u+d_v}\right)\left(\frac{\sqrt{d_u+d_v}}{\sqrt{d_ud_v}}\right)\right)^2}{\sum_{uv \in E}\frac{\sqrt{d^2_u+d^2_v}}{d_u+d_v}}\\
&= \frac{\left(\sum_{uv \in E}\left(\frac{1}{\sqrt{d_ud_v}}\right)\left(\frac{\sqrt{(d_u^2+d_v^2)(d_u+d_v)}}{d_u+d_v}\right)\right)^2}{DSO(G)}\\
&\geq \frac{\left(\sum_{uv \in E}\left(\frac{1}{\sqrt{d_ud_v}}\right)\left(\frac{\sqrt{(d_u+d_v)^2}}{d_u+d_v}\right)\right)^2}{DSO(G)}\\
&=\frac{(R(G))^2}{DSO(G)}.
\end{align*}
\qed

\begin{theorem}\label{theorem15}
Let $G$ be a connected graph of size $m$ with the maximum degree $\Delta$ and the minimum degree $\delta\geq 2$. Then
$$\frac{\sqrt{2}}{2\Delta}SDD(G)\leq DSO(G)\leq \left(\frac{\Delta}{2\delta}\right)\sqrt{mSDD(G)}.$$
Equality holds if and only if $G$ is regular graph. 
\end{theorem}
\proof
By substituting $r=1$, $x_i=\frac{\sqrt{d^2_u+d^2_v}}{d_u+d_v}$ and $a_i=\frac{d_ud_v}{(d_u+d_v)^2}$ into Lemma \ref{lemma6}, we have
\begin{align*}
SDD(G)&=\sum_{uv \in E}\frac{d_u^2+d_v^2}{d_ud_v}=\sum_{uv \in E} \left(\frac{\sqrt{d_u^2+d_v^2}}{d_u+d_v}\right)^2\left(\frac{(d_u+d_v)^2}{d_ud_v}\right)\\
&=\sum_{uv \in E}\frac{\left(\frac{\sqrt{d_u^2+d_v^2}}{d_u+d_v}\right)^2}{\frac{d_ud_v}{(d_u+d_v)^2}}\geq \frac{\left(\sum_{uv \in E}\frac{\sqrt{d_u^2+d_v^2}}{d_u+d_v}\right)^2}{\sum_{uv \in E}\frac{d_ud_v}{(d_u+d_v)^2}}\\
&=\frac{(DSO(G))^2}{\sum_{uv \in E}\frac{d_ud_v}{(d_u+d_v)^2}}
\end{align*}
Since for any $u \in V$, $\delta \leq d_u \leq \Delta$, we have $\frac{\delta^2}{4\Delta^2}\leq \frac{d_ud_v}{(d_u+d_v)^2}\leq \frac{\Delta^2}{4\delta^2}$. The equality holds if and only if $\delta=\Delta$. Therefore we get
\begin{align*}
SDD(G)&\geq \frac{(DSO(G))^2}{\sum_{uv \in E}\frac{d_ud_v}{d(d_u+d_v)^2}}\geq \frac{(DSO(G))^2}{\sum_{uv \in E}\frac{\Delta^2}{4\delta^2}}= \frac{(DSO(G))^2}{\left(\frac{\Delta}{2\delta}\right)^2m}.
\end{align*}
Consequently,
$$DSO(G)\leq \sqrt{\left(\frac{\Delta}{2\delta}\right)^2m SDD(G)}=\left(\frac{\Delta}{2\delta}\right)\sqrt{mSDD(G)}.$$
The equality holds if and only if $\delta=\Delta$ and hence, $G$ is a regular graph.\\
For the lower bound, by substituting $a_i=\sqrt{d^2_u+d^2_v}=b_i$ and $w_i=\frac{1}{d_u+d_v}$ into Lemma \ref{lemma7}, we get
\[
\left( \sum_{uv \in E} \frac{d_u^2+d_v^2}{d_u+d_v} \right) \left( \sum_{uv \in E} \frac{d_u^2+d_v^2}{d_u+d_v} \right)
\leq 
\max \left( 
\sqrt{2}\Delta \sum_{uv \in E}\frac{\sqrt{d_u^2+d_v^2}}{d_u+d_v},\,
\sqrt{2}\Delta \sum_{uv \in E}\frac{\sqrt{d_u^2+d_v^2}}{d_u+d_v}
\right)
\sum_{uv \in E} \frac{d_u^2+d_v^2}{d_u+d_v}.
\]
Hence, 
$$\sum_{uv \in E} \frac{d_u^2+d_v^2}{d_u+d_v}  \leq \sqrt{2}\Delta \sum_{uv \in E}\frac{\sqrt{d_u^2+d_v^2}}{d_u+d_v}.$$
Therefore, since $d_u\geq 2$ for any $u\in V$
\begin{align*}
\sqrt{2}\Delta DSO(G) &= \sqrt{2}\Delta \sum_{uv \in E}\frac{\sqrt{d_u^2+d_v^2}}{d_u+d_v}\geq  \sum_{uv \in E} \frac{d_u^2+d_v^2}{d_u+d_v}\\
&\geq \sum_{uv \in E} \frac{d_u^2+d_v^2}{d_ud_v}=SDD(G),
\end{align*}
and this gives the desired inequality. The equality holds if and only if $d_u=d_v$ for any $uv \in E$. Consequently, the equality holds if and only if $G$ is a regular graph. 
\qed

\begin{theorem}\label{theorem16}
Let $G$ be a connected graph with the maximum degree $\Delta$ and the minimum degree $\delta$. Then
$$\frac{SO(G)}{2\Delta}\leq DSO(G) \leq \frac{SO(G)}{2\delta}.$$
Equality holds if and only if $G$ is regular graph.
\end{theorem}
\proof
Since for any $uv \in E$, $2\delta \leq d_u+d_v \leq 2\Delta$, then
$$\frac{1}{2\Delta}\sum_{uv \in E}\sqrt{d_u^2+d_v^2}\leq \sum_{uv \in E}\sqrt{d_u^2+d_v^2}\times \frac{1}{d_u+d_v}\leq \frac{1}{2\delta}\sum_{uv \in E}\sqrt{d_u^2+d_v^2}.$$
$$\Longrightarrow \frac{SO(G)}{2\Delta}\leq DSO(G) \leq \frac{SO(G)}{2\delta}.$$
The equality holds if and only if $\delta=\Delta$, that is, $G$ is a regular graph. 
\qed

\begin{theorem}\label{theorem19}
Let $G$ be a simple graph of size $m$ with the maximum degree $\Delta$. Then
$$DSO(G)\geq \frac{1}{2\Delta}\sqrt{F(G)+m(m-1)\left(\Pi_{F}(G)\right)^{\frac{1}{m}}}.$$
\end{theorem}
\proof
For $a_i=d_u^2+d_v^2$ in Lemma \ref{lemma8}, we get
\begin{align*}
\left(2\Delta DSO(G)\right)^2&\geq \left(\sum_{uv \in E}\left(d_u+d_v\right)\frac{\sqrt{d_u^2+d_v^2}}{d_u+d_v}\right)^2\\
&=\left(\sum_{uv \in E}\sqrt{d_u^2+d_v^2}\right)^2\\
&\geq \sum_{uv \in E}(d_u^2+d_v^2)+m(m-1)\left(\prod_{uv \in E}(d_u^2+d_v^2)\right)^{\frac{1}{m}}\\
&=F(G)+m(m-1)\left(\Pi_{F}(G)\right)^{\frac{1}{m}}.
\end{align*}
Therefore, we have
$$DSO(G)\geq \frac{1}{2\Delta}\sqrt{F(G)+m(m-1)\left(\Pi_{F}(G)\right)^{\frac{1}{m}}}.$$
\qed

\begin{theorem}
Let $G$ be a graph of order $n$ and size $m$. Then
$$\frac{1}{2}BSO(G) \leq DSO(G) \leq \frac{n}{4}BSO(G).$$
Equality on the left-hand side holds if and only if $G\simeq mK_2$ and in the right-hand side holds if and only if $G$ consists of components,
each of which is a regular graph (not necessarily of equal degree).
\end{theorem}
\proof
For each $uv \in E$, using Jensen's inequality with considering the convex function $f(x)=\frac{1}{x}$, we have
\begin{equation}\label{e1}
\frac{2}{d_u+d_v}\leq \frac{1}{2}\left(\frac{1}{d_u}+\frac{1}{d_v}\right)=\frac{1}{2}\left(\frac{d_u+d_v}{d_ud_v}\right).
\end{equation}
Therefore, using (\ref{e1}) and since $d_u+d_v \leq n$ for any $uv \in E$, we get
\begin{align*}
DSO(G)&=\sum_{uv \in E}\frac{\sqrt{d_u^2+d_v^2}}{d_u+d_v}=\frac{1}{2}\sum_{uv \in E}\left(\sqrt{d_u^2+d_v^2}\right)\left(\frac{2}{d_u+d_v}\right)\\
&\leq \frac{1}{4}\sum_{uv \in E}\left(\sqrt{d_u^2+d_v^2}\right)\left(\frac{d_u+d_v}{d_ud_v}\right)=\frac{1}{4}\sum_{uv \in E}\left(d_u+d_v\right)\left(\frac{\sqrt{d_u^2+d_v^2}}{d_ud_v}\right)\\
&\leq \frac{n}{4}\sum_{uv \in E}\frac{\sqrt{d_u^2+d_v^2}}{d_ud_v}=\frac{n}{4}BSO(G).
\end{align*}
Equality holds if and only if $d_u=d_v$ for any $uv \in E$ in the graph $G$. \\
Since $0<\frac{1}{d_u}\leq 1$ for any $u \in V$, thus we have $\frac{1}{d_u}+\frac{1}{d_v}\leq 2$. Therefore, we get $\frac{2}{d_u+d_v}\geq \frac{1}{d_ud_v}$. 
\begin{align*}
DSO(G)&=\sum_{uv \in E}\frac{\sqrt{d_u^2+d_v^2}}{d_u+d_v}=\frac{1}{2}\sum_{uv \in E}\left(\sqrt{d_u^2+d_v^2}\right)\left(\frac{2}{d_u+d_v}\right)\\
&\geq\frac{1}{2}\sum_{uv \in E}\frac{\sqrt{d_u^2+d_v^2}}{d_ud_v}=\frac{1}{2}BSO(G).
\end{align*}
Equality holds if and only if $d_u=1$  for any $u \in V$. 
\qed

\vspace*{0.5 cm}
\section{Conclusion}
In this paper, we have conducted a comprehensive investigation of the Diminished Sombor index (DSO) and its relationships with several classical topological indices. By establishing new bounds and characterizing extremal graphs, our results contribute to a deeper understanding of the mathematical properties and potential applications of the DSO index, particularly in chemical graph theory. The connections drawn between DSO and established indices such as the Zagreb, Albertson, Harmonic, Randi\'c, and geometric-arithmetic indices highlight the versatility and relevance of this new invariant.

Looking ahead, further research could focus on extending the theoretical framework of the DSO index to more general graph classes, exploring its computational aspects, and investigating its practical performance in chemical and interdisciplinary applications. Such studies may reveal new insights and broaden the utility of the DSO index in both mathematical and applied contexts.

\noindent\textbf{Acknowledgements} The present study was supported by Golestan University, Gorgan, Iran (research number: 1749). The author truly appreciates Golestan University for this support.

\end{document}